\documentclass[11pt]{amsart}

\usepackage{amssymb}
\def\r{{\mathbb R}}
\def\n{{\mathbb N}}
\def\z{{\mathbb Z}}
\def\c{{\mathbb C}}

\def\t{\times}
\def\d{\delta}
\def\rr{{\mathcal R}}

\begin{document}

\centerline{On $p$-adic Speh representations}
\centerline{A. I. Badulescu}
\ \\

\ \ \ \ \ \ \ \ \ \ \ \ \ \ \ \ \ \ \ \ \ \ \ \ \ \ \ \ \ \ \ \ \ \ \ \ \ \ \ \ \ \ \ \ \ \ \ \ \ \ \ \ \ \ \ \ \ \ \ \ \ \ \ \ \ \ \ \ \ \ \ \ \ \ \ \ \ \ \ \ \ \ \ \ \ \ \ \ \ \ \ \ \ \ \ \ \ \ \ \ {\footnotesize{\it A la Virgen del Roc\'io}}\\
\ \\
\ \\
\ \\
{\bf Abstract:} {\it This note contains simple proofs of some known results {\rm (}unitarity, character formula{\rm )} on Speh representations of a group $GL_n(D)$ where $D$ is a local non Archimedean division algebra of any characteristic.}\\
\ \\
\ \\
{\bf Introduction.} In this note I give simple proofs of some known results on Speh representations of groups $GL_n(D)$ where $D$ is a central division algebra of finite dimension over a local non archimedean field $F$ of any characteristic. The new idea is to use the Moeglin-Waldspurger algorithm (MWA) for computing the dual of an irreducible representation. For the unitarizability of Speh representations, previous proofs were based either on the trace formula and the close fields theory, or on deep results of type theory. The proof here is "combinatoric", independent of $D$, the characteristic, and Bernstein's (also called U0) theorem. Short general proofs using MWA are also given for other known facts. The proof is always the same: one wants to prove a relation (R) involving some representations (for example an induced representation $\pi$ is irreducible). One starts by writing the "naive" relation (R') between these representations, known from standard theory, but not as strong as (R) (for example the semi simplification of $\pi$ is a sum with non negative coefficients of some irreducible representations). Usually (R') has more terms that (R), because it is weaker, and one wants to prove that some terms which are supposed to appear in (R') are actually not there (for example all the subquotients of $\pi$ except of the expected one have coefficient zero). The method then is to consider also the dual relation (R'') to (R'), and to play with the MWA, in order to show that the mild constraints one has on (R') and (R'') are enough to show that the extra terms are not there and (R') reduces actually to (R).

All the results in this paper are already known, and here I only give new short proofs. So I kindly ask the reader, when using one of these facts, to quote, at least in first place, the original reference (see the historical notice at the end of the paper).
Beside the Zelevinsky and Tadi\'c classification of the admissible dual ([Ze], [Ta2]), the proofs here rely on [Au] (dual of an irreducible representation is irreducible), [MW] and [BR2] (algorithm for computing the dual) and some easy tricks from [Ta1], [Ba3] and [CR], and do not involve any complicated technique.

The idea of searching for "simple proof" for classification of unitary representations, together with a list of basic tricks to use, are due to Marko Tadi\'c ([Ta3] for example). He also was the first to formulate some properties of Speh representations  ({\it formula for ends of complimentary series, character formula, Speh representations are prime elements in the ring of representations, dual of a Speh representation is Speh}) and to prove them when $D=F$. The starting point of my proof here of the assertion {\it Speh representations are unitary} is also due to Tadi\'c who found the simple but brilliant trick reducing the problem of unitarity to a problem of irreducibility.

I would like to thank Guy Henniart who read the paper and made useful observations.\\
\ \\
{\bf Notation.}
Let $(F,|\ |_F)$ be a local non archimedean field of any characteristic $ch(F)$ and $D$ a central division algebra of finite dimension $d^2$ over $F$. Let $\nu_n :GL_n(D)\to \c^\t$ be the character $g\mapsto |\det(g)|_F$, denoted simply $\nu$ when no confusion may occur ($\det$ is the reduced norm). All representations here will be admissible of finite length.  If $\pi$, $\pi'$ are representations of $GL_a(D)$ and $GL_b(D)$ such that $a+b=n$, let $P$ be the parabolic subgroup of $GL_n(D)$ containing the upper triangular matrices and of Levi factor  ${GL_a(D)}\t GL_b(D)$ and set $\pi\t\pi':=ind_P^{GL_n(D)}\pi\otimes \pi'$. 

The results in this section for which no other reference is given may be found in [Ze] (for $D=F$) and [Ta2] (for general $D$). In [Ta2] the characteristic zero is assumed, but this restriction may be removed ([BHLS] 2.2).
Let $Gr_n$ denote the Grothedieck group of admissible representations of finite length of $GL_n(D)$, and, by convention, $Gr_0=\z$. Set $\rr:=\oplus_{n\in\n} Gr_n$. The composition law $\t$ induces a composition law $*$ in $\rr$ which makes $\rr(+,*)$ into a commutative $\z$-algebra ($\t$ is not commutative, but $*$ is). We write {\it representation} for {\it isomorphy class of representations} where no confusion may occur. The set $Irr$ of irreducible representations of $GL_n(D)$, $n\in\n^\t$, is a natural linear basis of $\rr$.

Let  $\d$ be an essentially square integrable (i.e. twist of a square integrable by a character) representation of $GL_n(D)$. Let $s(\d)$ be the smallest positive real number $x$ such that $\d\times \nu^x\d$ is reducible. Then $s(\d)\in\n^*$ and $s(\d)|d$ (so $s(\d)=1$ if $D=F$). Set $\nu_\d:=\nu^{s(\d)}$. Then $\delta$ may be obtained as the unique irreducible subrepresentation of a representation of type $\nu_\d^{a}\rho\ \t\  \nu_\d^{a-1}\rho\ \t \ \nu_\d^{a-2}\rho\ \t\  ...\ \t\  \nu_\d^{a-m+1}\rho$, with $a\in\r$, $m$ a divisor of $n$, and $\rho$ a cuspidal unitary representation of $GL_{\frac{n}{m}}(D)$. Also, $a,m$ and $\rho$ are determined by $\d$. We then write $\d=Z(\rho,a,m)$. We have $s(\rho)=s(\d)$. We set $e(\d):=a-\frac{m-1}{2}$. Then $\d$ is unitary (i.e. is square integrable) if and only if $e(\d)=0$.
In this case we simply write $\d=Z(\rho,m)$, understood as: $\d$ is the unique subrepresentation of $\nu_\d^{\frac{m-1}{2}}\rho\ \t\  \nu_\d^{\frac{m-1}{2}-1}\rho\ \t \ \nu_\d^{\frac{m-1}{2}-2}\rho\ \t\  ...\ \t\  \nu_\d^{-\frac{m-1}{2}}\rho$. We make the convention that $Z(\rho,a,0)$ is the unit element of $\rr$.

To any irreducible representation $\pi$ of $GL_n(D)$ we may associate a unique multiset $[\d_1,\d_2,...,\d_k]$ such that 

- $\d_i$ is an essentially square integrable representations of some $GL_{n_i}(D)$, $n_i\geq 1$,

- $\sum_i n_i=n$ and 

- if we reorder the $d_i$ such that $e(\d_i)$ is decreasing with $i$, then $\pi$ is the unique irreducible quotient of $\d_1\t \d_2\t ...\t \d_k$.

We write $\pi=L([\d_1,\d_2,...,\d_k])$. 
We call {\it standard representation} an element of $\rr$ which is a product $\d_1*\d_2*...*\d_k$ with $d_i$ essentially square integrable representations. It is known that $\d_1*\d_2*...*\d_k$ determines the multiset $[\d_1,\d_2,...,\d_k]$, so we may write $\pi=L(\d_1*\d_2*...*\d_k)$ instead of $\pi=L([\d_1,\d_2,...,\d_k])$. Let ${\mathcal D}$ be the set of essentially square integrable representations in $\rr$ and $Std$ be the set of standard representations. Then $L$ realizes a bijection from $Std$ to $Irr$ (and its restriction to ${\mathcal D}$ is the identity).
If $\Theta$ is a standard representation and $\pi=L(\Theta)$, we then say that $\pi$ is {\it the Langlands quotient} of $\Theta$.\\  

If $\d$ is square integrable and $k\in \n^*$, then $u(\d,k)$ (resp. $\bar{u}(\d,k)$) will denote the Langlands quotient of \ $\nu_\d^{\frac{k-1}{2}}\d\ \t\  \nu_\d^{\frac{k-3}{2}}\d\ \t \ \nu_\d^{\frac{k-5}{2}}\d\ \t\  ...\ \t\  \nu_\d^{-\frac{k-1}{2}}\d$ (resp. of \ $\nu^{\frac{k-1}{2}}\d\ \t\  \nu^{\frac{k-3}{2}}\d\ \t \ \nu^{\frac{k-5}{2}}\d\ \t\  ...\ \t\  \nu^{-\frac{k-1}{2}}\d$). The representations $u(\d,k)$ are called {\bf Speh representations}.

If $\d=Z(\rho, l)$, then $\d_+$ (and $\d_-$) will stand for $Z(\rho,l+1)$ (and $Z(\rho,l-1)$).\\

Let $\iota$ be the Aubert involution on $\rr$ from [Au] (see also [ScS]) which we will call here {\it duality}. It takes an irreducible representation of $GL_n(D)$ to plus or minus an irreducible representation of $GL_n(D)$. We will denote $i$ the involution of the set of irreducible representations induced by $\iota$ (i.e. forgetting the sign; for $D=F$ it is the Zelevinsky involution from [Ze]). We recall that $\iota (\pi)*\iota (\pi')=\iota (\pi*\pi')$.\\

We call {\it segment} a subset of $\r$ of the form $S=\{a,a+1,a+2,...,a+m-1\}$, $m\in\n$, with the convention $S=\emptyset$ if $m=0$. Set $l(S)=m$ and call it the {\it length of $S$}. We say $a$ is the {\it beginning} of $S$ and $a+m-1$ is the {\it ending} of $S$. 

We call {\it multisegment} a non empty multiset $M=[S_1,S_2,...,S_k]$ of non empty segments. 
Set $m(M):=\max_{S\in M}l(S)$ (the {\it maxlength} of $M$) and $t(M):=k$ (the {\it thickness} of $M$). Let $B(M)$ (resp. $E(M)$) be the multiset of the beginnings (resp. endings) of segments in $M$. We call {\it support} of $M$ the multiset reunion of all the segments in $M$. Two segments  $S_1$ and $S_2$ are said to be {\it linked} if  $S_1\not\subseteq S_2$, $S_2\not\subseteq S_1$ and $S_1\cup S_2$ is a segment. If $M$ and $M'$ are multisegments we say $M'$ is obtained from $M$ by an {\it elementary operation} if there exists a pair of linked segments $S_1, S_2$ in $M$ such that $M'$ is obtained from $M$ by replacing $[S_1,S_2]$ with $[S_1\cup S_2, S_1\cap S_2]$ where it is understood we forget $S_1\cap S_2$ in case it is empty. One introduces then a partial order on multisegments by putting $M\leq M'$ if $M$ is obtained from $M'$ through a finite number of elementary operations. It is easy to see that on the partially ordered set of multisegments $m$ is decreasing, while $t$, $B$ and $E$ are increasing (for the relation on sets $X\leq Y$ if $X\subset Y$).

Let $R$ be the set of unitary cuspidal representations in $\rr$. If $\rho\in R$, let $K_{\rho}$ be the set $\{\nu_\rho^{\alpha} \rho,\ \alpha\in\r\}$. Let $Irr_\rho$ ($Std_\rho$, ${\mathcal D}_\rho$) be the subset of $Irr$ ($Std$, $\mathcal D$) made of representations with cuspidal support included in $K_\rho$. If $\rho,\rho'\in R$ and $\rho\neq \rho'$, if $\pi\in Irr_\rho$ and $\pi'\in Irr_{\rho'}$, then $\pi * \pi'\in Irr$ ([Ze], [Ta2]). Let $\rr_\rho$ be the linear span of $Irr_\rho$ in $\rr$. Then $\rr_\rho$ is a subalgebra of $\rr$ and one has $\rr=\t_{\rho\in R}\rr_\rho$. Each  $\rr_\rho$ is stable by the involution $\iota$.

Let us fix $\rho\in R$ and for the rest of this paragraph consider only representations in $\rr_\rho$. Then essentially square integrable representations  are parametrized by segments and so standard  representations  are parametrized by multisegments. The (partial) order relation defined on the set of multisegments induces an order relation on $Std_\rho$. The bijection $L: Std \to Irr$ induces a bijection from $Std_\rho$ to $Irr_\rho$ and we get by transfer a partial order on $Irr_\rho$. We will use this bijection $L$ to compare for example an irreducible representation to a standard representation. Also we will use $E(\pi)$, $t(\pi)$ etc.  for an irreducible representation $\pi$, being understood that this is $E(M)$, $t(M) $ etc., where $M$ is the multisegment associated to the standard representation $\Theta$ such that $\pi=L(\Theta)$.  For any standard representation $\Theta$, one has in $\rr$ (actually in $\rr_\rho$)
$$\Theta=L(\Theta)+\sum_{\pi_j<L(\Theta)}a_i\pi_j,$$ 
$a_i\in\n^*$, $\pi_j\in Irr$ (this is an avatar of the Langlands quotient theory explicitly stated in this stronger form in [Ze] and [Ta2] Theorem 5.3). 

As a consequence, $Std_\rho$ is a linear basis of the sub-space $\rr_\rho$.
So $(\rr_\rho,+,*)$ is isomorphic with the ring of polynomials over $\z$ with (commutative) variables in the set ${\mathcal D}_\rho$. A standard representation $\Theta\in \rr_\rho$ is then a monomial. Notice that its degree equals $t(L(\Theta))$ and is also $\geq t(\pi_j)$ for all the $\pi_j<L(\Theta)$ (because $t$ is decreasing). Decomposing an irreducible representation on the standard basis comes to write it as a polynomial.\\
\ \\
{\bf Important fact.} If $\pi\in \rr_\rho$ has a decomposition $\pi=\sum_i m_i\pi_i$ with $m_i\neq 0$ and $\pi_i$ distinct irreducible representations, then, as polynomials, we have $\deg(\pi)=\max_i \deg(\pi_i)$. To see this one writes $\pi_i=L(\Theta_i)$ and notice that the degree of $\pi_i$ is also the degree of the monomial $\Theta_i$. If, among the $\pi_i$ {\it of maximal degree as polynomials}, $\pi_{i_0}$ is maximal for the order relation on representations, then the monomial $\Theta_{i_0}$ appears with non zero coefficient in $\pi$, which shows that $\deg(\pi)\geq \max_i \deg(\pi_i)$. The inequality $\leq$ is obvious.\\

Using the relation $\rr=\t_{\rho\in R}\rr_\rho$ one may see that $Std$  is a linear basis of $\rr$ and $\rr$ is a polynomial algebra with variables in ${\mathcal D}$. Every element of $Std$ may be written in a unique way as a product of elements of $Std_\rho$, $\rho\in R$, which are almost all trivial. So we get a product order relation on $Std$. This order relation passes then to $Irr$. If $\Theta\in Std$, one has then in $\rr$ 
$$\Theta=L(\Theta)+\sum_{\pi_j<L(\Theta)}a_i\pi_j$$ 
with $a_i\in\n^*$ and $\pi_j\in Irr$. And if $\Theta=\prod_{\rho\in A} \Theta_\rho$ where $A$ is a finite subset of $R$ and $\Theta_\rho\in Std_\rho$,then $L(\Theta)=\prod_{\rho\in A} L(\Theta_\rho)$. 

Let $\pi,\pi'\in Irr$ and set $\pi=L(\Theta)$, $\pi'=L(\Theta')$ with $\Theta, \Theta'\in Std$. Then 
$$\pi *\pi' =x+\sum_{\pi_j< x}a_i\pi_j,$$ 
where $x=L(\Theta*\Theta')$, $a_i\in\n$, $\pi_j\in Irr$. We will call $x$ {\it the Langlands quotient} of $\pi*\pi'$.\\ 
\ \\
{\bf The algorithm.} For this section we fix $\rho\in R$ and we parametrize $Irr_\rho$ by multisegments. 
Let $\pi\in Irr_\rho$. If $M$ and $M^\#$ are the multisegments of $\pi$ and $i(\pi)$, Moeglin and Waldspurger gave an algorithm for computing $M^\#$ from $M$ ([MW] for $D=F$, [BR2] for $D\neq F$). I will call here MWA this algorithm. Here is a description:

If $\d=\{b,b+1,...,e\}$ is a non empty segment, we set $\d^-=\{b,b+1,...,e-1\}$ (which may be empty).

Let $M$ be a multisegment. We associate to $M$ the multisegment
$M^{\#}$ in  the following way~: let $e$ be the greatest number in $E(M)$. Chose:

- a shortest segment $\d_{0}\in M$ with ending $e$.  

- a shortest segment $\d_{1}\in M$ with ending $e-1$ and not included in $\d_{0}$

- a shortest segment $\d_{2}\in M$ with ending $e-2$ and not included in $\d_{1}$ etc..\\
Stop when it is not possible to find the next segment. Say $\d_r$ is the last one. Then $S_0=\{e-r,e-r+1,...,e\}$ is the first
segment of $M^{\#}$. Let now
$M^-$ be the multisegment obtained from $M$ after replacing $\d_i$ with $\d_i^-$ for all $i\in\{0,1,...,r\}$ (if one of these is empty, just erase it).
Starting with $M^-$ what we
have  done with $M$, we find the second segment of $M^{\#}$, and so on
(so in the end we have that $M^{\#}$ is the multiset union of
$[S_0]$  and $(M^-)^{\#}$). As observed in [MW], if the greatest number in $E(M)$ is $e$, then {\it all} the segments of $M^\#$  ending with $e$ will be constructed exclusively with points from $E(M)$. As an easy consequence, we have the following property:

{\bf (P)} {\it If $E(M)$ is included in some segment $S$ of length $k$ {\rm (}i.e. $\forall x\in E(M), x\in S${\rm )}, then $l(M^\#)\leq k$}. (Indeed, all the segments of $M^\#$ ending with the greatest point $e$ of $S$ will be included in $S$, so will have length $\leq k$. And once we are done with $e$, the set of endings will be (as a consequence) included in the shift of $S$ with $-1$.)

Also, every segment of $M^\#$ being constructed taking elements on distinct segments of $M$, we obviously have 

{\bf (P')} : {\it $l(M^\#))\leq t(M)$, i.e. the maxlength of $M^\#$ is less than or equal to the thickness of $M$}.\\
\ \\
{\bf The results.} In the sequel we will consider elements $\pi$ of $Irr$ which will actually be, most of the time, obviously in some $Irr_\rho$. For these, we will loosely say things like "let $M$ be the multisegment of $\pi$", or apply MWA, being understood we parametrized by segments. 
The assertions {\bf (b), (c), (d)} in the following theorem are the {\bf (U1), (U2), (U3)} of [Ta2].\\
\ \\
{\bf Theorem 1.} Let $\d$ be a square integrable representation of some $GL_n(D)$.

{\it {\bf (a)} If $\d=Z(\rho,l)$ then $i(u(\d,k))=u(\tau,l)$ where $\tau:=Z(\rho,k)$.

{\bf (b)} The representation $u(\d,k)$ is unitary.

{\bf (c)} The representations $\pi(u(\d,k),\alpha):=\nu_\d^\alpha u(\d,k)\times \nu_\d^{-\alpha} u(\d,k)$ are irreducible and unitary for $\alpha\in ]-\frac{1}{2},\frac{1}{2}[$.

{\bf (d)} $u(\d,k)$ is a prime element of the ring $\rr$.}\\
\ \\
{\bf Proof} of {\bf (a):} is obvious by MWA.\\
\ \\
{\bf Proof} of {\bf (b):} Tadi\'c ([Ta1]) reduced it by induction to the proof of :

\ \ \ \ \ \ \ \ \ \ \ \ {\bf (i)} for all $k\geq 1$, $u(\d,k)\t u(\d,k)$ is irreducible and

\ \ \ \ \ \ \ \ \ \ \ \ {\bf (ii)} for all $k\geq 2$, $u(\d,k-1)\t u(\d,k+1)$ is irreducible.
Here I prove {\bf (i)} and {\bf (ii)} (this trick of Tadi\'c was also used in [BHLS] 4.1):

{\bf Proof} of {\bf (i):} let $\d=Z(\rho,l)$ and set $\tau:=Z(\rho,k)$. Write

$$u(\d,k)* u(\d,k)=\pi+\sum_{\pi_j<\pi} a_j\pi_j$$
where $\pi$ is the Langlands quotient of the product and $a_j\in\n$. The same way:

$$u(\tau,l)* u(\tau,l)=\pi'+\sum_{\pi'_j<\pi'} b_j\pi'_j.$$ 
By {\bf (a)}, $u(\tau,l)$ is the dual of $u(\d,k)$. The second relation is then dual to the first. The MWA implies that the Langlands quotient $\pi'$ is the dual of $\pi$. If some $a_j\neq 0$, then $\pi_j$ appears effectively in the first sum and its dual, which we call $\pi'_j$, appears in the second sum by linear independence of the characters. As $l(\tau)=k$ and $\pi'_j<\pi'$, we have $m(\pi'_j)>k$ (at the first elementary operation the longest segment will already be longer than $k$; then $m$ is decreasing). This is impossible by {\bf (P)}, since $E(\pi_j)\subseteq E(u(\d,k)* u(\d,k)$  ($E$ is increasing) which is included in a segment of length $k$. So $a_j=0$ for all $j$.

{\bf Proof} of {\bf (ii):} write

$$u(\d,k-1)* u(\d,k+1)=\pi+\sum_{\pi_j<\pi} a_j\pi_j$$
where $\pi$ is the Langlands quotient of the product and $a_j\in\n$. Then write the dual relation (using {\bf (a)} for the left side) :

$$u(\tau_-,l)* u(\tau_+,l)=\pi'+\sum_{\pi'_j<\pi'} b_j\pi'_j$$ 
where MWA implies that the Langlands quotient $\pi'$ is the dual of $\pi$. Let $\pi_j$, $\pi'_j$ be two dual representations appearing effectively in the sums. As $l(\d)=l$ we know $m(\pi_j)\geq l+1$. Now 
$$E(\pi'_j)\subset E(u(\tau_+,l)* u(\tau_-,l))=x+[0,1,1,2,2,3,3,...,l-1,l-1,l]$$
and
$$B(\pi'_j)\subset B(u(\tau_+,l)* u(\tau_-,l))=y+[0,1,1,2,2,3,3,...,l-1,l-1,l]$$
where for the sake of the reader I shifted with $x=\frac{k-l-1}{2}$ and $y=\frac{-k-l+1}{2}$. Notice that $x+l$ always belongs to $E(\pi'_j)$, so if $M$ is the multisegment of $\pi'_j$, then the construction of (the first segment of) $M^\#$ starts with $x+l$. Now I will show that {\it the first segment $S_0$ of $M^\#$ {\rm (}i.e. the one ending in $x+l${\rm )} does not contain $x$}. Indeed, the segment $\d_0$ of $\pi'_j$ containing $x+l$ cannot be shorter  than the one of $u(\tau_+,l)* u(\tau_-,l)$ containing $x+l$ (elementary operations cannot shorten {\it this one}). So the beginning of $\d_0$ is $\leq y+l-1$. But the sequence of beginnings of segments $\d_0,\d_1,\d_2,...$ used in the construction of $S_0$ has to be strictly decreasing, and there is no place in $B(\pi'_j)\backslash\{y+l\}$ for such a sequence with $l+1$ terms.
Hence the first segment of $M^\#$ is constructed only with points from $x+\{1,2,...,l\}$ so its length is $\leq l$. But, as this one does not contain $x$, we have that $E(M^-)$ is included (in the sense of {\bf (P)}) in the segment $x+\{0,1,2,...,l-1\}$, which implies by {\bf (P)} that the other segments of $M^\#$ are also of length $\leq l$. This is in contradiction with $l(\pi_j)\geq l+1$.\\
\ \\
{\bf Proof {\rm of} (c):} Tadi\'c showed in [Ta2] irreducibility for $\alpha\in ]-\frac{1}{2},\frac{1}{2}[\backslash\{0\}$ and also that unitarity would follow from irreducibility in $\alpha=0$, which has been proved in {\bf (i)} (see the proof of {\bf (b)}).\\
\ \\
{\bf Proof {\rm of} (d):} Assume $u(\d,k)$ is {\it not} a prime element of $\rr$. Then it is not a prime element of $\rr_\rho$ (easy) and we will work in this ring. 
As $\rr_\rho$ is a polynomial ring over $\z$, it is a UFD, hence $u(\d,k)$ is reducible. Write 
$$u(\d,k)=S* T$$
where $S$ and $T$ are non scalar elements of $\rr_\rho$. Then
$$u(\tau,l)=\iota (S)* \iota (T).$$

As a polynomial, $u(\d,k)$ has degree $k$ and so  $deg(S),deg(T)<k$. Then, if we write $S$, $T$ as linear combination of irreducible representations, every irreducible representation $\pi$ in $S$ or $T$ has $t(\pi)<k$. So every irreducible representation $i(\pi)$ in $\iota (S)$ or $\iota (T)$ has $m(i(\pi))<k$ (by {\bf (P')}). But, if $\alpha$ is an irreducible representation in $\iota (S)$, then $t(\alpha)\leq deg(\iota (S))$ and if $\beta$ is an irreducible representation in $\iota (T)$ then $t(\beta)\leq deg(\iota (T))$ so finally $t(\alpha)+t(\beta)\leq l$ (because $deg(\iota (S))+deg(\iota (T))=deg(u(\tau,l))=l$). But, as all the segments of $\alpha$ and $\beta$ are strictly shorter than $k$, $\alpha*\beta$ lives on a group of smaller size than $u(\tau,l)$. It is impossible for $u(\tau,l)$ to be a sum of such representations.\qed
\ \\
\ \\
{\bf Corollary 1.} {\it Let ${\mathcal T}$ be the set of all finite products of {\rm (}irreducible unitary{\rm )} representations of type $u(\d,k)$ and $\pi(u(\d,k),\alpha)$ for $k\in\n^*$, $\d$ square integrable and $\alpha\in ]0,\frac{1}{2}[$. Let $\mathcal U$ be the set of all irreducible unitary representations of all $GL_m(D)$, $m\in\n^*$. Then ${\mathcal T}\subset {\mathcal U}$.}\\
\ \\
{\bf Proof.} By prop. 2.13 in [Ba3], a product $u_1\t u_2\t ...\t u_k$ of irreducible unitary representations is irreducible unitary if, for all $1\leq i\leq k$, $u_i\t u_i$ is irreducible. By {\bf (i)} (see the proof of {\bf (b)}), $U=u(\d,k)\t u(\d,k)$ is indeed irreducible. To show  that $\pi(u(\d,k),\alpha)\t \pi(u(\d,k),\alpha)$ is also irreducible we write it as $\nu_\d^\alpha U\t \nu_\d^{-\alpha} U$.  Now $\nu_\d^\alpha U$ and $\nu_\d^{-\alpha} U$ are irreducible representations, and for $0< \alpha < \frac{1}{2}$ no segment of the first one is linked to no segment of the second one. So $\nu_\d^\alpha U\t \nu_\d^{-\alpha} U$ is irreducible by Proposition 2.2 and Lemma 2.5 of [Ta2].\qed
\ \\

If $\d=Z(\rho,l)$ is a square integrable representation and $k\in \n^*$, set
$$F(\d,k):=\nu^{-\frac{k+l}{2}}(\sum_{w\in W_k^l}(-1)^{sgn(w)}\prod_{i=1}^kZ(\rho,i,w(i)+l-i))\in\rr,$$
where $W_k^l$ is the set of permutation $w$ of $\{1,2,...,k\}$ such that $w(i)+l\geq i$ for all $i\in\{1,2,...,k\}$.\\
\ \\
{\bf Theorem 2.} {\it

{\bf (e)} If $\d=Z(\rho,l)$ then one has, in $\rr$, \ $u(\d,k)=F(\d,k).$

{\bf (f)} One has \ \ $\nu_\d^{-\frac{1}{2}} u(\d,k) * \nu_\d^{\frac{1}{2}} u(\d,k)=u(\d,k-1)* u(\d,k+1) + u(\d_-,k) * u(\d_+,k).$}\\
\ \\
{\bf Proof {\rm of} (e):} the case $k=1$ is trivial, $k=2$ is known by [Ze] and [Ta2]. Assume by induction that $u(\d,k)=F(\d,k)$ and $u(\d,k-1)=F(\d,k-1)$ for any square integrable representation $\d$.

Let us write 
$$\nu_\d^{-\frac{1}{2}} u(\d,k)* \nu_\d^{\frac{1}{2}} u(\d,k)=u(\d,k-1)* u(\d,k+1) + u(\d_-,k) * u(\d_+,k)+\sum_{\pi_j<u(\d,k-1)* u(\d,k+1)} a_j\pi_j$$
and taking the dual:
$$\nu_\d^{-\frac{1}{2}} u(\tau,l)* \nu_\d^{\frac{1}{2}} u(\tau,l)=u(\tau,l-1)* u(\tau,l+1) + u(\tau_-,l) * u(\tau_+,l)+ \sum_{\pi'_j<u(\tau,l-1)* u(\tau,l+1)} b_j\pi'_j.$$
These formulas need some explanation: as $u(\d,k-1)* u(\d,k+1)$ is irreducible by {\bf (ii)} (see the proof of {\bf (b)}), it is the Langlands quotient of $\nu_\d^{-\frac{1}{2}} u(\d,k)* \nu_\d^{\frac{1}{2}} u(\d,k)$. For the same reason, $u(\tau,l-1)* u(\tau,l+1)$ is the Langlands quotient of $\nu_\d^{-\frac{1}{2}} u(\tau,l)* \nu_\d^{\frac{1}{2}} u(\tau,l)$.  Now, as $i(u(\tau,l-1)* u(\tau,l+1))=u(\d_-,k) * u(\d_+,k)$, this representation appears with multiplicity one in $\nu_\d^{-\frac{1}{2}} u(\d,k)* \nu_\d^{\frac{1}{2}} u(\d,k)$ by duality. As $i(u(\d,k-1)* u(\d,k+1))=u(\tau_-,l) * u(\tau_+,l)$, this representation appears with multiplicity one in 
$\nu_\d^{-\frac{1}{2}} u(\tau,l)* \nu_\d^{\frac{1}{2}} u(\tau,l)$ by duality. So the two $\Sigma$ are dual to one another and we also have the property {\bf [A]}: {\it the coefficient $b_j$ of  $\pi'_j=u(\tau_-,l) * u(\tau_+,l)$ is zero}.

Chenevier and Renard ([CR]) recognized in $F(\d,k)$ a determinant with entries in $\rr$ and proved using Lewis Carroll's identity that  
$${\bf [F]}\ \ \ \ \ \ \ \ \ \ \ \ \nu_\d^{-\frac{1}{2}} F(\d,k)* \nu_\d^{\frac{1}{2}} F(\d,k)=F(\d,k-1)* F(\d,k+1) + F(\d_-,k) * F(\d_+,k).$$

Using the induction assumption we replace to obtain:
$$\nu_\d^{-\frac{1}{2}} u(\d,k)* \nu_\d^{\frac{1}{2}} u(\d,k)=u(\d,k-1)* F(\d,k+1) + u(\d_-,k) * u(\d_+,k).$$
Subtracting from our first formula we get: 
$$\sum_{\pi_j<u(\d,k-1)* u(\d,k+1)} a_j\pi_j=u(\d,k-1)\big(F(\d,k+1)-u(\d,k+1)\big).$$ 

So taking the dual, $u(\tau_-,l)$ divides $\sum_{\pi'_j<u(\tau,l-1)* u(\tau,l+1)} b_j\pi'_j$. Write $\sum_{\pi'_j<u(\tau,l-1)* u(\tau,l+1)} b_j\pi'_j=u(\tau_-, l)* \sum c_p\theta_p$ with $c_p\in\z^*$ and $\theta_p$ distinct irreducible representations. Chose a maximal $\theta_p$. Then $L(u(\tau_-,l)* \theta_p)$ is one of the $\pi'_j$ such that $b_j\neq 0$. But $t(\pi'_j)\leq 2l$ so $t(\theta_p)\leq l$.  As $E(\pi_j)$ is included in a segment of length $k+1$, $m(\pi'_j)\leq k+1$ (by {\bf (P)}) and so $m(\theta_p)\leq k+1$. As the support of $\theta_p$ has $l(k+1)$ elements and we proved it contains less then $l$ segments, all of length less than $k+1$, the only possibility is $\theta_p=u(\tau_+,l)$. So $u(\tau_-,l)* \theta_p= u(\tau_-,l) * u(\tau_+,l)$. As this representation is irreducible, we have $L(u(\tau_-,l)* \theta_p)= u(\tau_-,l) * u(\tau_+,l)$.This is in contradiction with {\bf [A]}. So what we called $\sum c_p\theta_p$ is void, and this implies {\bf (e)} at level $k+1$.\\
\ \\
{\bf Proof {\rm of} (f):} The formula {\bf [F]} implies that {\bf (e)} for all $k$ implies {\bf (f)} ([CR]).\qed
\ \\
\ \\
{\bf Corollary 2.} {\it Let $\rm C$ be the Jacquet-Langlands transfer for square integrable representations from {\rm [DKV]} between $GL_{nd}(F)$ to $GL_n(D)$ and $\rm LJ$ the transfer of irreducible representations from {\rm [Ba2]} between $GL_{ndk}(F)$ to $GL_{nk}(D)$.  If $\d$ is a square integrable representation of $GL_{nd}(F)$, then ${\rm LJ}(u(\d,k))=\bar{u}({\rm C}(\d),k)$.}\\
\ \\
{\bf Proof.} It is proved  in [Ta5] to be an easy consequence of the character formula {\bf (e)}. Notice the [DKV] transfer is known also when $ch(F)\neq 0$ ([Ba1]).\qed

\ \\
{\bf Remark.} Bernstein's theorem asserts that {\it induced representations from unitary irreducible representations remain irreducible}. Bernstein proved this theorem for $D=F$ in [Be]. Sécherre proved it for general $D$ in [Se]. Tadi\'c showed ([Ta1]) using Bernstein's theorem that actually 
${\mathcal T}= {\mathcal U}$ (the Corollary 2 here is just one inclusion) and this is the only proof of the equality we have up to now. The proof of the Bernstein's theorem is difficult, especially for general $D$, and a new short one  would be most welcome. This theorem has not been used in this paper.\\ 
\ \\
{\bf Historical notice.}\\
\ \\
For $D=F$, Theorem 1 and Corollary 1 have been proved by Tadi\'c in [Ta1] (1986), {\bf (f)} in [Ta7] (1987), and {\bf (e)} in [Ta4] (1995). Tadi\'c used Bernstein's theorem ([Be]) and Zelevinsky derivatives. He conjectured this results should hold also for $D\neq F$ in [Ta2] (1990), where he proved {\bf (d)} in this case.\\ 
\ \\
For $D\neq F$, $ch(F)=0$:

- {\bf (b)} has been proved in [BR1] (2004) for $D$ of characteristic zero (as it requires a result from [Ba2] based on the Arthur trace formula).

- Tadi\'c gave then a proof of all the other results in this paper in [Ta5] (2006) assuming only Bernstein's theorem were true for $D\neq F$ (but not involving anymore Zelevinsky's derivatives which are specific to $D=F$). I gave myself a proof of these facts in [Ba3] (2008) without using the Bernstein's theorem which was not known for $D\neq F$ at that time. Eventually Sécherre proved Bernstein's theorem for $D\neq F$ in [Se] (2009).\\ 
\ \\
For proofs in the case $D\neq F$ and $ch(F)\neq 0$ : [BHLS] (2010).\\
\ \\
A proof of similar assertions for $D$ archimedean (i.e. $D=\r,\c,{\mathbb H}$) may be found in the literature. See [Sp], [Vo], [Ta6], and [BR3]. It would be long to include the statement here because the analogue of the Zelevinsky classification is different, so notation and formulation change.\\
\ \\
\ \\
{\bf Bibliography.}

[Au] A.-M. Aubert, Dualité dans le groupe de Grothendieck
de la catégorie des représentations lisses de longeur
finie d'un groupe réductif $p$-adique, {\it Trans. Amer. Math. Soc.},
347 (1995), 2179-2189 (there is an Erratum: {\it Trans. Amer. Math. Soc.},
348 (1996), 4687-4690).

[Ba1] A. I. Badulescu, Correspondance de Jacquet-Langlands en
caractéristique non nulle, {\it Ann. Scient. \'Ec. Norm. Sup.} 35
(2002), 695-747.

[Ba2] A. I. Badulescu, Jacquet-Langlands et unitarisabilité, {\it J. Inst. Math. Jussieu}
6(3):349-379, 2007.

[Ba3] A. I. Badulescu, Global Jacquet-Langlands correspondence, multiplicity one and classification
of automorphic representations {\it Invent. Math}., 172(2):383-438, 2008. With
an appendix by Neven Grbac.

[Be] J. N. Bernstein, $P$-invariant distributions on $GL(N)$ and the
classification of unitary representations of $GL(N)$
(non-Archimedean case), in {\it Lie groups and representations II},
Lecture Notes in Mathematics 1041, Springer-Verlag, 1983.

[BR1] A. I. Badulescu, D. Renard, Sur une conjecture de Tadi\'c,
{\it Glasnik Matematicki} 39 (59) no. 1 (2004), p. 49-54.

[BR2] A. I. Badulescu, D. Renard, Zelevinsky involution and
Moeglin-Waldspurger algorithm for $GL(n,D)$, in {\it Functional
analysis IX} (Dubrovnik, 2005), 9-15., Various Publ. Ser. 48, Univ.
Aarhus, Aarhus, 2007.

[BR3] A. I. Badulescu, D. Renard, Unitary dual of $GL_n$ at archimedean places and global Jacquet-Langlands correspondence.
{\it Compositio Math.} 146, vol. 5 (2010), p. 1115-1164. 

[BHLS] A. I. Badulescu, G. Henniart, B. Lemaire, V. Sécherre, Sur le dual unitaire de $GL(r,D)$.
{\it American Journal of Math.}, vol. 132, no. 5 (2010), p. 1365-1396.

[CR] G. Chenevier, D. Renard, Characters of Speh representations and Lewis Caroll identity, {\it  Represent. Theory}  12  (2008), 447--452. 

[DKV] P. Deligne, D. Kazhdan, and M.-F. Vign´eras. Représentations des algèbres centrales
simples p-adiques. In {\it Representations of reductive groups over a local field}, Travaux
en Cours, pages 33 - 117. Hermann, Paris, 1984.

[MW] C. Moeglin, J.-L.Waldspurger, Sur l'involution de Zelevinski,
{\it J. Reine Angew. Math.}  372  (1986), 136-177.

[ScS] P. Schneider, U. Stuhler, Representation theory and sheaves on
the Bruhat-Tits building,  {\it Inst. Hautes \'Etudes Sci. Publ.
Math.}  No. 85 (1997), 97-191.

[Se] V. S\'echerre, Proof of the Tadi\'c U0 conjecture on the
unitary dual of $GL(m,D)$, J. Reine Angew. Math.  626  (2009), 187-203.

[Sp] B. Speh, Unitary representations of $Gl(n, R)$ with nontrivial (g, K)-cohomology, {\it Invent.
Math.} 71(3):443-465, 1983.

[Ta1] M. Tadi\'c, Classification of unitary representations in
irreducible representations of a general linear group
(non-Archimedean case), {\it Ann. Scient. \'Ec. Norm. Sup.} 19
(1986), 335-382.

[Ta2] M. Tadi\'c, Induced representations of $GL(n;A)$ for a
$p$-adic division algebra $A$, {\it J. Reine angew. Math.} 405
(1990), 48-77.

[Ta3] M. Tadi\'c, An external approach to unitary representations,
{\it Bulletin Amer. Math. Soc.} 28, No. 2 (1993), 215-252.

[Ta4] M. Tadi\'c, On characters of irreducible unitary
representations of general linear groups, {\it Abh. Math. Sem. Univ.
Hamburg} 65 (1995), 341-363.

[Ta5] M. Tadi\'c, Representation theory of $GL(n)$ over a $p$-adic division algebra and unitarity
in the Jacquet-Langlands correspondence, {\it Pacific Journal of Math.} 223 (2006), no 1, 167-200.

[Ta6] M. Tadi\'c, $\widehat{GL(n,\c)}$ and $\widehat{GL(n,\r)}$,  {\it Automorphic forms and $L$-functions II. Local aspects},  285--313, Contemp. Math., 489, Amer. Math. Soc., Providence, RI, 2009.

[Ta7] M. Tadi\'c, Topology of unitary dual of non-Archimedean ${\rm GL}(n)$,  {\it Duke Math. J.}  55  (1987),  no. 2, 385-422.

[Vo] D. A. Vogan, Jr. The unitary dual of $GL(n)$ over an Archimedean field, {\it Invent. Math.}
83(3):449-505, 1986.

[Ze]  A. Zelevinsky, Induced representations of reductive $p$-adic
groups II, {\it Ann. Scient. \'Ec. Norm. Sup.} 13 (1980), 165-210.
\ \\
\ \\
Alexandru Ioan Badulescu\\
Université de Montpellier 2,\\
Département de mathématiques, case courrier 051,\\
Place Eugène Bataillon, 34095 Montpellir Cedex, FRANCE\\
E-mail: badulesc@math.univ-poitiers.fr

\end{document}